\documentclass[10pt]{book}
\usepackage[sectionbib]{natbib}
\usepackage{array,epsfig,fancyheadings,rotating}
\usepackage{algorithm, algorithmic}

\usepackage{amsmath}
\usepackage{amssymb}
\usepackage{amsfonts}
\usepackage{multirow}
\usepackage{amsthm}
\usepackage{bm}

\newtheorem{theorem}{Theorem}
\newtheorem{lemma}{Lemma}
\newtheorem{corollary}{Corollary}
\newtheorem{proposition}{Proposition}
\newtheorem{condition}{Condition}
\theoremstyle{definition}

\newcommand{\e}{\bm{\varepsilon}} 
\newcommand{\E}{\mathbb{E}} 
\newcommand{\argmin}{\mathop{\rm argmin}\limits}
\newcommand{\prob}{\mathbb{P}}
\begin{document}


\renewcommand{\baselinestretch}{1.2}

\markright{ \hbox{\footnotesize\rm Submitted to Statistica Sinica
}\hfill\\[-13pt]
\hbox{\footnotesize\rm
}\hfill }

\markboth{\hfill{\footnotesize\rm SHOTA KATAYAMA AND HIRONORI FUJISAWA} \hfill}
{\hfill {\footnotesize\rm SPARSE AND ROBUST LINEAR REGRESSION} \hfill}

\renewcommand{\thefootnote}{}
$\ $\par


\fontsize{10.95}{14pt plus.8pt minus .6pt}\selectfont
\vspace{0.8pc}
\centerline{\large\bf SPARSE AND ROBUST LINEAR REGRESSION:}
\vspace{2pt}
\centerline{\large\bf AN OPTIMIZATION ALGORITHM AND}
\vspace{2pt}
\centerline{\large\bf ITS STATISTICAL PROPERTIES}
\vspace{.4cm}
\centerline{Shota Katayama and Hironori Fujisawa}
\vspace{.4cm}
\centerline{\it Tokyo Institute of Technology and The Institute of Statistical Mathematics}
\vspace{.55cm}
\fontsize{9}{11.5pt plus.8pt minus .6pt}\selectfont


\begin{quotation}
\noindent {\it Abstract:}
This paper studies sparse linear regression analysis with outliers in the responses.
A parameter vector for modeling outliers is added to the standard linear regression model and then 
the sparse estimation problem for both coefficients and outliers is considered.
The $\ell_{1}$ penalty is imposed for the coefficients, while various penalties including redescending type penalties 
are for the outliers.
To solve the sparse estimation problem, we introduce an optimization algorithm.
Under some conditions, we show
the algorithmic and statistical convergence property for the coefficients obtained by the algorithm.
Moreover, it is shown that the algorithm can recover the true support of the coefficients with probability going to one.
\par

\vspace{9pt}
\par
\end{quotation}\par

\def\thefigure{\arabic{figure}}
\def\thetable{\arabic{table}}

\fontsize{10.95}{14pt plus.8pt minus .6pt}\selectfont

\setcounter{chapter}{1}
\setcounter{equation}{0} 
\noindent {\bf 1. Introduction}

Linear regression with a large number of covariates is a general and fundamental problem in recent data analysis.
A standard method to overcome this problem is the 
least absolute shrinkage and selection operator (Lasso) proposed by \cite{tibshirani96}.
For a large number of covariates, 
it is natural to assume the sparsity which means that many of the covariates are not relevant to the responses.
The Lasso can draw relevant covariates automatically and simultaneously estimate the remaining coefficients to be zero.
From the theoretical point of view, the Lasso has been largely studied.
For instance, \cite{bickel09}, \cite{meinshausen09} and \cite{wainwright09}
gave the convergence rate under several norms and 
showed that the Lasso estimates can recover the true support of the coefficients.
See also \cite{efron2004}, \cite{zhao06}, \cite{zou06} and \cite{geer09}.

Recent linear regression analysis requires some complex structures in addition to a large number of covariates.
One of them is an outlier structure. This structure often appears in many applications such as
signal detection, image and speech processing, communication network and so on.
It is well known that the standard method which uses the $\ell_{2}$ loss
outputs inaccurate estimates when outliers exist.
The most popular way for robustifying against the outliers is to use the M-estimation procedure 
which replaces the $\ell_{2}$ loss by an other loss with a bounded influence function; e.g., 
the Huber's, the skipped-mean and the Hampel's robust loss (see, for instance, \cite{huber09} for more details).
However, optimizing such a robust loss function with a sparse penalty requires much computational cost.
The $\ell_{2}$ loss function is easier to deal with.
Recently, \cite{she11} proposed a novel approach for robust parameter estimation, using
an outlier model where a parameter vector for modeling outliers
is added to the standard linear regression model.
The corresponding $\ell_{2}$ loss function with a sparse penalty for the outlier parameter vector is optimized.
The connection between sparse penalties and robust loss functions was also illustrated.

\lhead[\footnotesize\thepage\fancyplain{}\leftmark]{}\rhead[]{\fancyplain{}\rightmark\footnotesize\thepage}

This paper studies sparse and robust linear regression based on their outlier model
from a theoretical point of view, which was not treated in \cite{she11}.
Some theoretical analyses were discussed by \cite{nguyen13}, but
only the $\ell_{1}$ penalty was used for the outlier parameters. 
Their results were derived from the fact that the resulting estimate is a global optimum.
As \cite{she11} showed, many penalties having good robustness are non-convex and then
the resulting estimate is often a local optimum.
A general theory including non-convex penalties is thus of great interest.

Main contributions in this paper are the following.
We consider a larger class of penalties for outlier parameters including non-convex penalties and then
derive some statistical properties.
To avoid the problem of local optima, we directly analyze estimated coefficients which an optimization algorithm outputs.
We provide the upper bound of its $\ell_{2}$ error, which is divided by an algorithmic error and a statistical error. 
It is also shown that the algorithm can recover the true support of the coefficients.
Thus, this paper bridges a gap between the statistical theory and the computation algorithm,
which arises in using non-convex penalties.

The remainder of this paper is organized as follows.
We introduce the model and the optimization algorithm in Section 2.
In Section 3, theoretical analyses for the estimated coefficients which the algorithm outputs are provided.
We report numerical performances in Section 4. All the proofs are in Section 5.

Throughout the paper, 
a bold symbol denotes a matrix or a vector and 
its element is written by the fine symbol, e.g., $\bm{A}=(a_{ij})$ for $\bm{A}\in\mathbb{R}^{m\times n}$ and $\bm{a}=(a_{1},\dots, a_{m})^{T}$ for $\bm{a}\in\mathbb{R}^{m}$.
For any vector $\bm{a} \in \mathbb{R}^{m}$ and $1\le q < \infty$, we define the $\ell_{q}$ norms as 
$\|\bm{a}\|_{q}=(\sum_{i=1}^{m} |a_{i}|^{q})^{1/q}$ and define the $\ell_{0}$ and $\ell_{\infty}$ norms as
$\|\bm{a}\|_{0}=|\{i|\,a_{i}\neq 0\}|$ and $\|\bm{a}\|_{\infty} = \max_{1\le i\le m}|a_{i}|$, respectively.
Given a set $S \subset \{1,\dots, m\}$, 
$\bm{a}_{S}$ denotes the sub-vector having the elements of $\bm{a}$ 
corresponding to the set $S$, that is,
$\bm{a}_{S}=\{a_{i}|\;i\in S\}$.
For any two vectors $\bm{a}, \bm{b} \in \mathbb{R}^{m}$, $\langle\bm{a},\bm{b}\rangle$ denotes
the standard inner product.
We define the matrix $\ell_{1}$ norm as
$\|\bm{A}\|_{\ell_{1}} = \max_{1\le j\le n}\sum_{i=1}^{m}|a_{ij}|$ and for any symmetric matrix
$\bm{B}\in\mathbb{R}^{m\times m}$, we define its largest eigenvalue as $\xi_{max}(\bm{B})$.
For two positive sequences $a_{n}, b_{n}$ depending on $n$, 
the notation $a_{n}=O(b_{n})$ means that there exists a finite constant $C>0$ such that
$a_{n} \le Cb_{n}$ for a sufficiently large $n$, while $a_{n}=\Omega(b_{n})$ means that
$a_{n} \ge Cb_{n}$. Also, the notation $a_{n}=o(b_{n})$ means that $a_{n}/b_{n} \to 0$ as $n$ goes to infinity.

\setcounter{chapter}{2}
\setcounter{equation}{0} 
\vskip 14pt
\noindent {\bf 2. Sparse and robust linear regression}

\vskip 14pt
\noindent {\bf 2.1. Model and parameter estimation}

Consider the linear regression model with outliers
\begin{align}\label{2.1}
\bm{y}=\bm{X}\bm{\beta^{*}} + \sqrt{n}\bm{\gamma^{*}} + \bm{\varepsilon},
\end{align}
where $\bm{y}=(y_{1},\dots, y_{n})^{T}$ is an $n$ dimensional response vector, 
$\bm{X}=(x_{ij}) = (\bm{x}_{1},\dots, \bm{x}_{n})^{T}$ is an $n\times p$ covariate matrix,
$\bm{\beta}^{*}$ is a $p$ dimensional unknown coefficient vector,
$\bm{\gamma}^{*}$ is an $n$ dimensional unknown vector whose nonzero elements correspond to outliers
and $\bm{\varepsilon}$ is an $n$ dimensional random error vector.
In the model (\ref{2.1}), we assume that the $\ell_{2}$ norm of columns of $\bm{X}$ is $\sqrt{n}$.
Correspondingly, the coefficient of $\bm{\gamma}^{*}$ is assumed to be $\sqrt{n}$ 
to match its scale with the columns of $\bm{X}$.
The model (\ref{2.1}) can be found also in \cite{she11} and \cite{nguyen13}.
Our purpose is to estimate $\bm{\beta}^{*}$ accurately even when the number of covariates is large.
In this case, it is natural to assume that $\bm{\beta}^{*}$ has many zero elements (sparse).
Moreover, we can assume that $\bm{\gamma}^{*}$ is also sparse since the number of outliers is usually not large.
For this sparse structure, we introduce sparse penalties for both
coefficients and outliers. The parameters are estimated by solving the optimization problem
\begin{align}\label{2.2}
\argmin_{\bm{\beta}\in\mathbb{R}^{p}, \bm{\gamma}\in\mathbb{R}^{n}}
\frac{1}{2n}\|\bm{y} - \bm{X}\bm{\beta} - \sqrt{n}\bm{\gamma}\|_{2} + \lambda_{\beta}\sum_{j=1}^{p}w_{\beta, j}|\beta_{j}|
+\lambda_{\gamma}\sum_{i=1}^{n}w_{\gamma, i}P(\gamma_{i}),
\end{align}
where $\lambda_{\beta} > 0$ and $\lambda_{\gamma} > 0$ are tuning parameters 
for $\bm{\beta}$ and $\bm{\gamma}$, respectively and $P(\cdot)$ is a penalty function that encourages sparsity.
We often use a redescending type of $P(\cdot)$, since it can yield a small bias against strong outliers.
We consider the adaptive Lasso (\cite{zou06}) type optimization problem.
In (\ref{2.2}), $w_{\beta, i}$ and $w_{\gamma, j}$ are known weights.
Suppose that we have the preliminary estimators $\tilde{\bm{\beta}}=(\tilde{\beta}_{1},\dots, \tilde{\beta}_{p})^{T}$ and
$\tilde{\bm{\gamma}}=(\tilde{\gamma}_{1},\dots, \tilde{\gamma}_{n})^{T}$.
The weights are typically defined as $w_{\beta,j}=1/|\tilde{\beta_{j}}|$ and $w_{\gamma, i}=1/|\tilde{\gamma}_{i}|$.
For the details of the preliminary estimators used in this paper, see Section 3.3.

\vskip 14pt
\noindent {\bf 2.2. Optimization algorithm}

\begin{algorithm}
\caption{}
\label{algorithm1}
\begin{algorithmic}
\STATE {\bf Step 1.} Initialize $k \leftarrow 0$, $\bm{\beta}^{k} \leftarrow \bm{\beta}^{init}$ 
and $\bm{\gamma}^{k} \leftarrow {\rm argmin}_{\bm{\gamma}}L(\bm{\beta}^{0},\bm{\gamma})$.

\STATE {\bf Step 2.} Update $k \leftarrow k+1$,
\begin{align}
\bm{\beta}^{k} &\leftarrow {\rm argmin}_{\bm{\beta}} L(\bm{\beta}, \bm{\gamma}^{k-1}), \label{2.3}\\
\bm{\gamma}^{k} &\leftarrow {\rm argmin}_{\bm{\gamma}} L(\bm{\beta}^{k}, \bm{\gamma}). \label{2.4}
\end{align}

\STATE {\bf Step 3.} 
If they converge, then output current $(\bm{\beta}^{k}, \bm{\gamma}^{k})$ and stop the algorithm, otherwise return to 
Step 2.
\end{algorithmic}
\end{algorithm}

Let $L(\bm{\beta}, \bm{w})$ be the objective function in (\ref{2.2}).
To solve the optimization problem (\ref{2.2}), we introduce Algorithm \ref{algorithm1}.
This algorithm surely converges since it satisfies that 
\begin{align*}
L(\bm{\beta}^{0}, \bm{\gamma}^{0}) \ge L(\bm{\beta}^{1}, \bm{\gamma}^{0}) \ge L(\bm{\beta}^{1}, \bm{\gamma}^{1})\ge 
\cdots \ge 0
\end{align*}
from its construction.
We also note that the optimization problem in (\ref{2.3}) can be solved by a typical way
such as the coordinate descent algorithm (\cite{friedman10}).
The optimization problem in (\ref{2.4}) can be rewritten as 
\begin{align*}
\argmin_{\bm{\gamma}\in\mathbb{R}^{n}}\frac{1}{2n}\sum_{i=1}^{n}\big\{(y_{i} - \bm{x}_{i}^{T}\bm{\beta}^{k} - \sqrt{n}\gamma_{i})^{2} + \lambda_{\gamma}w_{\gamma, i}P(\gamma_{i})\big\}.
\end{align*}
Suppose that the problem ${\rm argmin}_{x}(z-x)^{2}/2 + \lambda P(x)$ has the explicit solution
$\Theta(z;\lambda)$. Then, (\ref{2.4}) can be written as
\begin{align*}
\gamma_{i}^{k} \leftarrow \frac{1}{\sqrt{n}}\Theta (y_{i} - \bm{x}_{i}^{T}\bm{\beta}^{k}; \lambda_{\gamma}w_{\gamma, i}),
\quad i=1,\dots, n,
\end{align*}
where $\gamma_{i}^{k}$ is the $i$th element of $\bm{\gamma}^{k}$.
Let this expression be denoted by $\bm{\gamma}^{k} \leftarrow h(\bm{\beta}^{k})$.
Then, the step 2 can be expressed with only the update of $\bm{\beta}$ as
\begin{align*}
\bm{\beta}^{k} \leftarrow \argmin_{\bm{\beta}\in\mathbb{R}^{p}}
L\big\{\bm{\beta}, h(\bm{\beta}^{k-1})\big\}.
\end{align*}

The function $\Theta(z;\lambda)$ is often called as the thresholding function. 
As seen in \cite{antoniadis01}, many sparse penalties, including the $\ell_{1}$, $\ell_{0}$ and
smoothly clipped absolute deviation (SCAD; \citet{fan97}) penalties, have an explicit solution $\Theta(z;\lambda)$. 
The $\ell_{1}$ penalty leads to the soft thresholding function $\Theta (z;\lambda) = {\rm sgn}(z)\max(|z| - \lambda, 0)$
where ${\rm sgn}(\cdot)$ denotes the sign function, and the $\ell_{0}$ penalty leads to the
hard thresholding function $\Theta(z;\lambda) = zI(|z| > \lambda)$, where $I(A)$ denotes the indicator function on the 
event $A$. For the SCAD penalty, the corresponding thresholding function is given by
\begin{align*}
\Theta (z;\lambda) = \begin{cases}
{\rm sgn}(z)(|z| - \lambda) & {\rm if\ }|z| \le 2\lambda \\
\frac{(a-1)z - a\lambda{\rm sgn}(z)}{a-2} & {\rm if\ }2\lambda < |z| \le a\lambda \\
z & {\rm if\ }|z| > a\lambda,
\end{cases}
\end{align*}
where $a=3.7$ is recommended by \citet{fan01}.

\vskip 14pt
\noindent {\bf 2.3. Connection to robust $M$-estimators}

Similarly in \citet{she11}, our procedure also has a close connection to robust $M$-estimators.
\begin{proposition}\label{prop2.1}
Let $\hat{\bm{\beta}}$ be the coefficient output of Algorithm \ref{algorithm1} and
let $\psi(z; \lambda) = z - \Theta (z;\lambda)$, where $\Theta (\cdot ;\lambda)$ is a thresholding function 
yielded from a penalty function $P(\cdot)$.
Then, for $j=1,\dots, p$, it follows that
\begin{align}\label{2.5}
\frac{1}{n}\sum_{i=1}^{n}x_{ij}\psi(y_{i} - \bm{x}_{i}^{T}\hat{\bm{\beta}}; \lambda_{\gamma}w_{\gamma, i})
+\lambda_{\beta}w_{\beta, j}\partial|\hat{\beta}_{j}| = 0,
\end{align}
where $\partial |\hat{\beta}_{j}|$ is the subgradient of $|\hat{\beta}_{j}|$ which means that
$\partial |\hat{\beta}_{j}|={\rm sgn}(\hat{\beta}_{j})$ if $\hat{\beta}_{j}\neq 0$, 
$\partial |\hat{\beta}_{j}| \in [-1,1]$ otherwise.
\end{proposition}

Proposition \ref{prop2.1} shows that the output $\hat{\bm{\beta}}$ satisfies 
the estimating equations (\ref{2.5}) with the $\psi$ function and the $\ell_{1}$ penalty.
Thus, our algorithm is closely related to the optimization problem
\begin{align}\label{2.6}
\argmin_{\bm{\beta}\in\mathbb{R}^{p}}\frac{1}{2n}\sum_{i=1}^{n}\Psi(y_{i} - \bm{x}_{i}^{T}\bm{\beta};\lambda_{\gamma}w_{\gamma,i}) +
\lambda_{\beta}\sum_{j=1}^{p}w_{\beta,j}|\beta_{j}|,
\end{align}
where $\frac{d}{dt}\Psi(t;\lambda)=\psi(t;\lambda)$.
It may take much computational cost to directly solve (\ref{2.6}),
particularly when the function $\Psi(\cdot; \lambda)$ is non-convex.
For the fast computation, we can use the first-order approximation of $\Psi(\cdot;\lambda)$,
but it loses the information of the original function.
For this reason, we use Algorithm \ref{algorithm1} instead of solving (\ref{2.6}).

The relationship between sparse penalties and robust loss functions are stated through 
the equation $\psi(z;\lambda) = z - \Theta(z;\lambda)$.
For instance, the $\ell_{1}$ penalty corresponds to the Huber's loss, the $\ell_{0}$ penalty to the
skipped-mean loss and the SCAD penalty to the Hampel's loss.
\citet{she11} gave their illustrations and more details.
One measure to characterize the robust loss function is the redescending property. It means that
$\psi(z;\lambda)$ goes to $0$, equivalently $|z-\Theta(z;\lambda)|$ goes to $0$, as $|z|$ goes infinity.
The soft thresholding function which corresponds to the $\ell_{1}$ penalty does not have 
the redescending property, but the hard and SCAD thresholding functions 
which correspond to the non-convex penalties have it.
Some other non-convex penalties including the non-negative garrote penalty (Garrote; \citet{gao98}) and 
the minimax concave penalty (MCP; \citet{zhang10}) also lead to that property.
Thus, penalties yielding good robustness are non-convex.

\vskip 14pt
\setcounter{chapter}{3}
\setcounter{equation}{0} 
\noindent {\bf 3. Theoretical Analysis}

The optimization problem (\ref{2.2}) with a non-convex penalty $P(\cdot)$ suffers from local minima.
We may be able to give a theoretical analysis for a global minimum of (\ref{2.2}), but
it is unclear which optimization algorithm can actually output the global minimum.
To avoid such a problem, we directly analyze the output of Algorithm \ref{algorithm1}.
Some theoretical analyses for computable solutions in linear regression model without outliers were provided by
\citet{zhang12}, \citet{fan13} and \citet{fan14}, but our analysis is essentially different from theirs.
They first derived some properties for a global minimum and then showed that
a computable solution shared those properties under additional conditions for the solution.

\vskip 14pt
\noindent {\bf 3.1. Notations}

First, we prepare some notations.
Let $S^{*} = {\rm supp}(\bm{\beta}^{*})=\{i|\,\beta_{i}^{*}\neq 0\}\subset \{1,\dots, p\}$
and $G^{*} = {\rm supp}(\bm{\gamma}^{*}) \subset \{1,\dots, n\}$.
These support sizes are written by $s^{*}=|S^{*}|$ and $g^{*}=|G^{*}|$.
For the preliminary estimators,
similarly let $\tilde{S}={\rm supp}(\tilde{\bm{\beta}})$ and $\tilde{G}={\rm supp}(\tilde{\bm{\gamma}})$
with sizes $\tilde{s}=|\tilde{S}|$ and $\tilde{g}=|\tilde{G}|$.
We define the restricted smallest eigenvalue $\delta_{min}$ as
\begin{align}\label{3.1}
\delta_{min}(u)=
\inf_{\|\bm{\delta}\|_{0}\le u}\frac{\|\bm{X}\bm{\delta}\|_{2}^{2}}{n\|\bm{\delta}\|_{2}^{2}} 
\end{align}
and the (doubly) restricted largest eigenvalue $\delta_{max}$ as
\begin{align}\label{3.2}
\delta_{max}(u,u')=
\sup_{\|\bm{\delta}\|_{0}\le u}\sup_{|G|\le u'}\frac{\|\bm{X}_{(G)}\bm{\delta}\|_{2}^{2}}{n\|\bm{\delta}\|_{2}^{2}}
\end{align}
where $\bm{X}_{(G)}$ is the sub-matrix having the rows of $\bm{X}$ corresponding to the set $G$, that is,
$\bm{X}_{(G)} = \{x_{ij}|\; i\in G, 1\le j\le p\}$.
The restricted eigenvalue is popular in the analysis of the Lasso (see, e.g., \citet{bickel09}),
but our $\delta_{max}$ is slightly different from the existing one.
The restriction is imposed on rows of $\bm{X}$ in addition to columns, corresponding to the outliers.
We shall provide an asymptotic analysis as $n$ goes to infinity.
In our theory, $p$, $s^{*}$ and $g^{*}$ may go to infinity depending on $n$.
Notice that the restricted eigenvalues also depend on $n$.

\vskip 14pt
\noindent {\bf 3.2. Properties of Algorithm \ref{algorithm1} output}

To derive a convergence property of Algorithm \ref{algorithm1},
we require the following conditions for the random error vector 
$\bm{\varepsilon} = (\varepsilon_{1},\dots, \varepsilon_{n})^{T}$, 
the thresholding function $\Theta(\cdot;\lambda)$ and the preliminary estimator 
$(\tilde{\bm{\beta}}^{}, \tilde{\bm{\gamma}}^{})$.
\begin{condition}\label{condition1}
{\rm 
The errors $\varepsilon_{1},\dots, \varepsilon_{n}$ are independently and identically distributed as
the zero mean sub-Gaussian distribution with a parameter $\sigma > 0$, that is, 
 $\E(\varepsilon_{i})=0$ and
$\E\{\exp(t\varepsilon_{i})\} \le \exp(t^{2}\sigma^{2}/2)$ for all $t\in\mathbb{R}$.
}
\end{condition}

\begin{condition}\label{condition2}
{\rm 
The thresholding function $\Theta(\cdot;\lambda)$ satisfies that
$\Theta(x;\lambda) = 0$ if $|x| \le \lambda$ and $|\Theta(x;\lambda) - x| \le \lambda$ for all $x\in\mathbb{R}$.
}
\end{condition}


\begin{condition}\label{condition3}
{\rm 
There exist a sequence $a_{n,1}\to 0$ and a constant $\kappa > 0$ such that
$\|\tilde{\bm{\beta}} - \bm{\beta}^{*}\|_{2}^{} + \|\tilde{\bm{\gamma}} - \bm{\gamma}^{*}\|_{2}\le \tilde{C}a_{n,1}$
and $\delta_{min}(\tilde{s}) \ge \kappa$ with probability going to one, where $\tilde{C}$ is some positive constant.
}
\end{condition}

Condition \ref{condition1} assumes that the errors in the model (\ref{2.1}) belong to
the sub-Gaussian distribution family. As seen in \citet{vershynin12}, this family covers
the Gaussian, Bernoulli and any distributions with a bounded support.
Condition \ref{condition2} decides a class of thresholding functions (or equivalently penalties)
which can be handled by our theory. For instance, this class includes the soft, hard, non-negative garrote, 
SCAD and MCP thresholding functions.
The first half claim in Condition \ref{condition3} implies consistent preliminary estimators at the rate $a_{n,1}$.
Under this condition, we have 
$|\tilde{\beta}_{j}|\ge|\beta_{j}^{*}| - |\tilde{\beta}_{j} - \beta_{j}^{*}| \ge |\beta_{j}^{*}| - \tilde{C}a_{n,1}$ for 
$j\in S^{*}$ by the triangle inequality. 
Hence, if $\min_{j\in S^{*}}|\beta_{j}^{*}| > \tilde{C}a_{n,1}$, then $|\tilde{\beta}_{j}| > 0$ for $j\in S^{*}$.
Similarly, we can show that if $\min_{i\in G^{*}} |\gamma_{i}^{*}| > \tilde{C}a_{n,1}$ then $|\tilde{\gamma}_{i}| > 0$
for $i\in G^{*}$. Thus, Condition \ref{condition3} leads to the screening property that
$S^{*} \subset \tilde{S}$ and $G^{*}\subset \tilde{G}$ if 
\begin{align}\label{screening}
\min\big\{\min_{j\in S^{*}}|\beta_{j}^{*}|, \min_{i\in G^{*}}|\gamma_{i}^{*}|\big\} > \tilde{C}a_{n,1}.
\end{align}
The second half claim in Condition \ref{condition3} 
is well known as the sparse Riesz condition (\citet{zhang08}) if $\tilde{s}$ is non-random.
When we use a preliminary estimator which has a non-random upper bound $s_{u}$, 
it reduces to $\delta_{min}(s_{u}) \ge \kappa$.
In Section 3.3, we shall clarify the order of $a_{n,1}$ and the non-random upper bound of $\tilde{s}$ for some preliminary estimators.
For a technical reason, we re-define the weights $w_{\beta,j}$ and $w_{\gamma,i}$
using a constant $R_{w}>0$ by
\begin{align}\label{rr}
w_{\beta,j}=\max\bigg(\frac{1}{|\tilde{\beta}_{j}|},\frac{1}{R_{w}}\bigg),\quad
w_{\gamma,i}=\min\bigg(\frac{1}{|\tilde{\gamma}_{i}|}, R_{w}\bigg)
\end{align}
for $i\in\tilde{G}$ and $j\in\tilde{S}$.
From these restrictions, it follows that $\min_{j\in \tilde{S}}w_{\beta,j} \ge R_{w}^{-1}$ and
$\max_{i\in\tilde{G}}w_{\gamma,i} \le R_{w}$.
They are required to exclude the case where $w_{\beta,j} \to 0$ and
$w_{\gamma,i}\to \infty$ as $n$ goes to infinity in our theory.
Although these are the same as the originals if $R_{w}$ is large, 
it is convenient for deriving the following theories.
Moreover, the same result holds if different constants 
$R_{w,\beta}$ and $R_{w,\gamma}$ are used in (\ref{rr}).

Under these conditions, we can show the algorithmic and statistical convergence property of Algorithm \ref{algorithm1}.
\begin{theorem}\label{th3.1}
Assume Conditions 1--3 and (\ref{screening}).
Let 
\begin{align*}
\lambda_{\beta} \ge 2CR_{w}\sqrt{\frac{\sigma^{2}\log p}{n}},\quad
\lambda_{\gamma} \le \frac{Cn}{R_{w}\max_{j\in\tilde{S}}\sum_{i\in\tilde{G}^{}}|x_{ij}|}\sqrt{\frac{\sigma^{2}\log p}{n}}
\end{align*} 
for a constant $C>\sqrt{2}$ and $R_{w}>0$ used in (\ref{rr}).
Then,
for any iteration $k\ge1$ and any initial value $\bm{\beta}^{init}$
 in Algorithm \ref{algorithm1}, 
\begin{align}\label{3.3}
\|\bm{\beta}^{k} - \bm{\beta}^{*}\|_{2}^{} 
\le \rho_{}^{k}\|\bm{\beta}^{init} - \bm{\beta}^{*}\|_{2} + 2\kappa^{-1}\sqrt{s^{*}}\lambda_{\beta}
\big(R_{w}^{-1}+\max_{j\in S^{*}}w_{\beta,j}\big)\sum_{i=0}^{k-1}\rho_{}^{i},
\end{align}
with probability going to one, 
where $\rho_{}=2\kappa^{-1}\delta_{max}(\tilde{s},\tilde{g})$.
\end{theorem}

The first term of the right hand side of (\ref{3.3}) can be regarded as the algorithmic error. 
It represents the effect of an initial value $\bm{\beta}^{init}$ in Algorithm \ref{algorithm1}.
This shows that if $\rho < 1$ then the effect vanishes from the bound exponentially as 
the number of iterations increases. Since $|\xi_{max}(\bm{M})| \le \|\bm{M}\|_{\ell_{1}}$
for any symmetric matrix $\bm{M}$, we have
\begin{align}\label{3.5}
\delta_{max}(\tilde{s},\tilde{g}) \le \max_{1\le i\le n; 1\le j\le p}x_{ij}^{2}\frac{\tilde{s}\tilde{g}}{n}.
\end{align}
Thus, if $\max_{ij}x_{ij}^{2}a_{n,2}^{2} = o(n)$, we have
$\rho < 1$ for sufficiently large $n$ with probability going to one, where
$a_{n,2}$ shall be defined later in Condition \ref{condition4}. 
The second term can be regarded as the statistical error.
We notice that $\max_{j\in S^{*}}w_{\beta,j} = O(1)$ with probability going to one
if $\min_{j \in S^{*}}|\beta_{j}^{*}| = \Omega(1)$.
In this case, the second term has 
$O(\sqrt{s^{*}}\lambda_{\beta})$, which is equivalent to the order
of the standard Lasso excluding the term $\sqrt{n}\bm{\gamma}^{*}$ from the model (\ref{2.1}) in advance.

Theorem \ref{th3.1} shows only the convergence result for the output.
Next, we shall show that the output can recover the true support.
We need an extra condition and a corollary.

\begin{condition}\label{condition4}
{\rm 
There exist a constant $C > 0$ and a sequence $a_{n,2}$, which may diverge, such that
$\max(\tilde{s},\tilde{g})\le C_{}a_{n,2}$ with probability going to one. 
}
\end{condition}

\begin{corollary}\label{cor3.1}
{\rm
Assume Conditions \ref{condition1}--\ref{condition4}, (\ref{screening}), $\max_{ij}|x_{ij}| = O(1)$,
$\min_{j\in S^{*}}|\beta_{j}^{*}|=\Omega(1)$ and $a_{n,2}^{2}=o(n)$. 
Let
$\lambda_{\beta}\ge C_{\beta}\{(\log p)/n)\}^{1/2}$ and
$\lambda_{\gamma}\le C_{\gamma}\frac{n}{a_{n,2}}\{(\log p)/n\}^{1/2}$ 
with some constants $C_{\beta} > 0$ and $C_{\gamma} > 0$.
If there exist some $k_{0} \ge 1$ and $C_{0} > 0$ such that 
$\prob(\|\bm{\beta}^{init}-\bm{\beta}^{*}\|_{2}\le C_{0}\rho^{-k_{0}}\sqrt{s^{*}}\lambda_{\beta}) \to 1$, then it follows that
$\prob(\|\bm{\beta}^{k} - \bm{\beta}^{*}\|_{2}\le C\sqrt{s^{*}}\lambda_{\beta})\to 1$ for any $k\ge k_{0}$
for some constant $C > 0$.
}
\end{corollary}

When $\max_{i,j}|x_{ij}|=O(1)$ and $a_{n,2}^{2}=o(n)$, we have $\rho = o(1)$. 
Consequently, Corollary \ref{cor3.1} can be immediately obtained from Theorem \ref{th3.1}.
The requirement 
$\prob(\|\bm{\beta}^{init}-\bm{\beta}^{*}\|_{2}\le C_{0}\rho^{-k_{0}}\sqrt{s^{*}}\lambda_{\beta}) \to 1$
is not so restrictive since $\rho = o(1)$. If we use $\tilde{\bm{\beta}}$ for the initial value $\bm{\beta}^{init}$, then
the condition $(a_{n,2}^{2}/n)^{k_{0}}a_{n,1}=O(\sqrt{s^{*}}\lambda_{\beta})$ is required.
From Corollary \ref{cor3.1}, we can show a theorem about the support recovery.
It should be noted that we only need to analyze the elements of $\bm{\beta}^{k}$ on $\tilde{S} \cap S^{*^{c}}$ 
since the screening property $S^{*} \subset \tilde{S}$ holds.

\begin{theorem}\label{th3.2}
Assume the same conditions as in Corollary \ref{cor3.1}. 
If
$\log p = O(\log n)$, 
$\sqrt{\log n}=o(\sqrt{n}\min_{i\in G^{*}}|\gamma_{i}^{*}|)$,
$a_{n,2}s^{*} = o(n)$
and $a_{n,1}\max\big(\sqrt{a_{n,2}s^{*}},g^{*}/\sqrt{n}\big) = o(1)$, then we have
$\prob\big(\bm{\beta}^{k}_{\tilde{S}\cap S^{*^c}} = \bm{0}\big) \to 1$ for any $k\ge k_{0}+1$, where
$k_{0}$ is defined in Corollary \ref{cor3.1}.
\end{theorem}

\vskip 14pt
\noindent {\bf 3.3. Preliminary estimator and its properties}

In the previous section, some statistical properties are derived for general preliminary estimators.
In this section, we introduce a concrete example and specify the orders
$a_{n,1}$ and $a_{n,2}$ in Conditions \ref{condition3} and \ref{condition4}.
Let us consider the Lasso type preliminary estimators of 
$\tilde{\bm{\theta}}=(\tilde{\bm{\beta}}^{T}, \tilde{\bm{\gamma}}^{T})^{T}$, given by
\begin{align}\label{pre}
\tilde{\bm{\theta}}=\argmin_{\bm{\theta}\in\mathbb{R}^{np}}\frac{1}{2n}\|\bm{y} - \bm{Z}\bm{\theta}\|_{2}^{2}
+\lambda_{\theta}\|\bm{\theta}\|_{1},
\end{align}
where $\lambda_{\theta} > 0$ is the tuning parameter
and $\bm{Z}=(\bm{X},\sqrt{n}\bm{I}_{n})$.
The estimator $\tilde{\bm{\theta}}$ is a simple variant of that proposed in \citet{nguyen13}.
Using different tuning parameters for $\bm{\beta}$ and $\bm{\gamma}$ may improve the accuracy of estimates.
However, even if the preliminary estimates do not have the high accuracy,
we can improve its accuracy in calculating (\ref{2.2}) with different tuning parameters.
For this reason, it would be enough to use the single tuning parameter in (\ref{pre}).

Based on the well-known analysis for the Lasso (see, e.g., \citet{bickel09} and \citet{wainwright09}), 
the following property can be shown.

\begin{proposition}\label{pr3.1}
Assume Condition \ref{condition1} and there exists a constant $\tilde{\kappa}>0$ such that
\begin{align}\label{3.8}
\min_{\bm{\theta}\neq \bm{0};\|\bm{\theta}_{U^{*^{c}}}\|_{1} \le 3\|\bm{\theta}_{U^{*}}\|_{1}}
\frac{\|\bm{Z}\bm{\theta}\|_{2}^{2}}{n\|\bm{\theta}\|_{2}^{2}} \ge \tilde{\kappa},
\end{align}
where $U^{*}={\rm supp}(\bm{\theta}^{*})$ with $\bm{\theta}^{*} = (\bm{\beta}^{*T},\bm{\gamma}^{*T})^{T}$.
Let $\lambda_{\theta} = C_{\theta}\{(\log \max(n,p))/n\}^{1/2}$ for sufficiently large $C_{\theta} > 0$.
Then, it follows that
\begin{gather}
\|\tilde{\bm{\theta}} - \bm{\theta}^{*}\|_{2} \le C\sqrt{\frac{(s^{*} + g^{*})\log\max(n,p)}{n}}, \label{3.9} \\
|{\rm supp}(\tilde{\bm{\theta}})| \le C\xi_{max}(\bm{Z}^{T}\bm{Z}/n)(s^{*} + g^{*}), \label{3.10}
\end{gather}
with probability going to one, where $C > 0$ is some constant.
\end{proposition}

Note that the condition (\ref{3.8}) is slightly different from that commonly used in the Lasso analysis.
It is involved in the extended matrix $\bm{Z}$ not in $\bm{X}$. 
For more details, see \citet{nguyen13}.
Remember Conditions \ref{condition3} and \ref{condition4}.
The bounds (\ref{3.9}) and (\ref{3.10}) will correspond to the orders $a_{n,1}$ and $a_{n,2}$, respectively.
However, the term $\xi_{max}(\bm{Z}^{T}\bm{Z}/n)$ may become troublesome since
it may diverge as $n$ or $p$ increases. 
To exclude it, after $\tilde{\bm{\theta}}$ is obtained with $\lambda_{\theta}$, 
we consider the threshold version of (\ref{pre})
with an additional tuning parameter $\tau_{\theta} > 0$ as
\begin{align}\label{thpre}
\tilde{\theta}_{j}^{th} = \tilde{\theta}_{j}I(|\tilde{\theta}_{j}| > \tau_{\theta}\lambda_{\theta}),\quad
j=1,\dots, n+p.
\end{align}
Now, let $\tilde{\bm{\theta}}^{th} = (\tilde{\theta}_{1}^{th},\dots, \tilde{\theta}_{n+p}^{th})^{T}$ and
$a_{n,1} = \{(s^{*}+g^{*})\log\max(n,p)/n)\}^{1/2}$.
Under the condition (\ref{screening}) in which $\tilde{C}$ is replaced by $2\tilde{C}$, 
we have $|\tilde{\theta}_{j}| > \tilde{C}a_{n,1}$ for any $j \in {\rm supp}(\bm{\theta}^{*})$.
Thus, if we select $\tau_{\theta}$ such that $\tau_{\theta}\lambda_{\theta} \le \tilde{C}a_{n,1}$, then
$\tilde{\theta}_{j}^{th} = \tilde{\theta}_{j}$ for any $j \in {\rm supp}(\bm{\theta}^{*})$.
Hence, the threshold version also holds the property (\ref{3.9}) with the same order as the original.
Meanwhile, note that
\begin{align*}
|{\rm supp}(\tilde{\bm{\theta}}^{th})\backslash {\rm supp}(\bm{\theta}^{*})|
&= \sum_{j\in {\rm supp}(\tilde{\bm{\theta}}^{th})\backslash {\rm supp}(\bm{\theta}^{*})}1
\le \sum_{j\not\in {\rm supp}(\bm{\theta}^{*})}
\frac{\tilde{\theta}_{j}^{2}}{\tau_{\theta}^{2}\lambda_{\theta}^{2}} \\
&\le\frac{\|\tilde{\bm{\theta}} - \bm{\theta}^{*}\|_{2}^{2}}{\tau_{\theta}^{2}\lambda_{\theta}^{2}}
\le C(s^{*} + g^{*}),
\end{align*}
if we select $\tau_{\theta} \ge C_{\tau}$ for sufficiently small $C_{\tau} > 0$,
which implies that $|{\rm supp}(\tilde{\bm{\theta}}^{th})| \le (1+C)(s^{*} + g^{*})$.
Here, the term $\xi_{max}(\bm{Z}^{T}\bm{Z}/n)$ is excluded.
Note that the conditions $\tau_{\theta}\lambda_{\theta}\le \tilde{C}a_{n,1}$ and $\tau_{\theta} \ge C_{\tau}$
are compatible since the order of $\lambda_{\theta}$ is smaller than or equal to that of $a_{n,1}$.

\vskip 14pt
\setcounter{chapter}{4}
\setcounter{equation}{0} 
\noindent {\bf 4. Numerical performance}

\begin{figure}[ht]
  \centering
  \includegraphics[width=1\textwidth, bb=0 0 799 432]{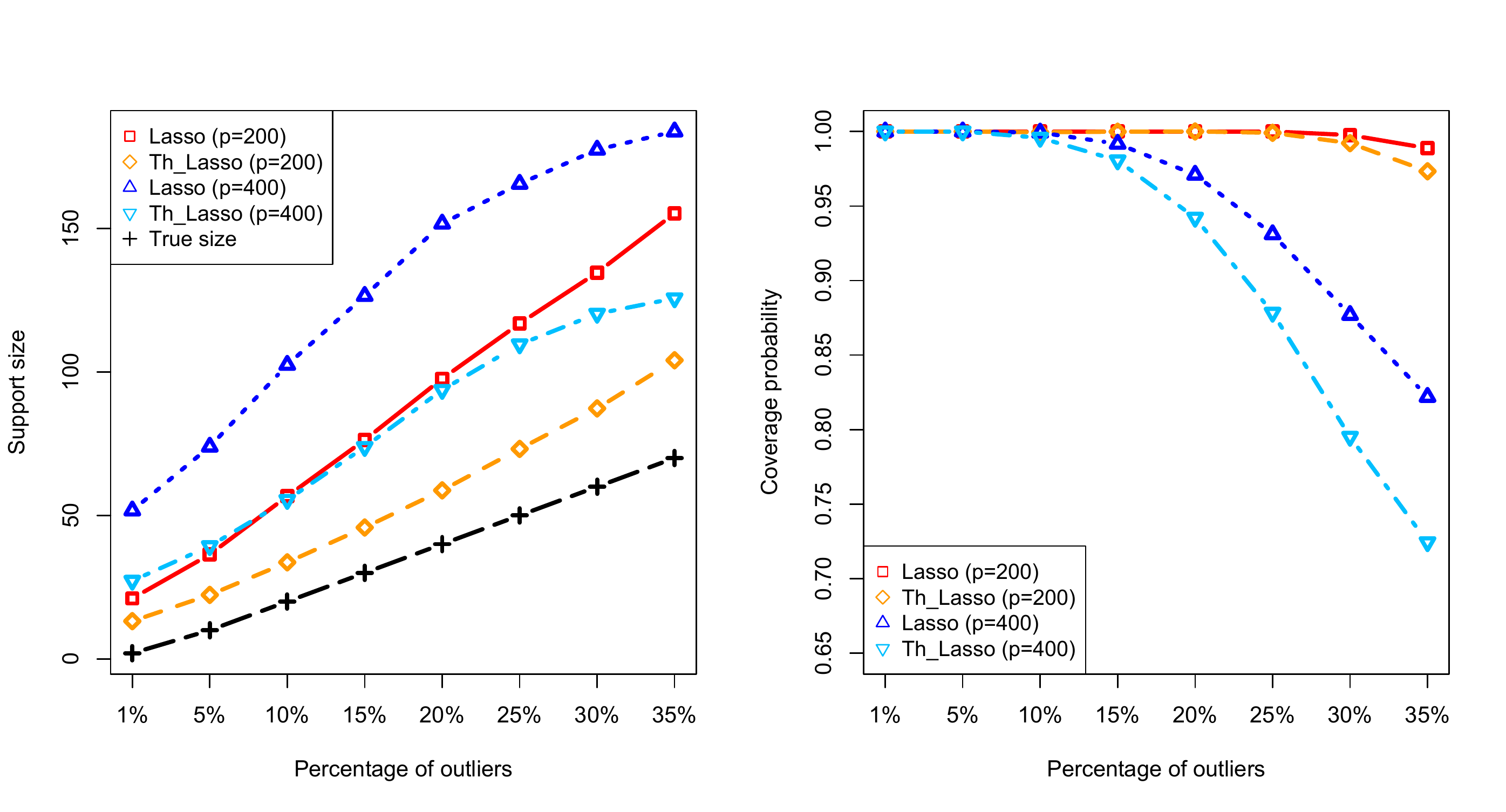}
  \caption{The support size (left) and the coverage probability (right) for the two preliminary estimators.
  ``Lasso" means (\ref{pre}) and ``Th\_Lasso" means (\ref{thpre}).
  Each point on the curve shows the mean based on 100 Monte Carlo simulations.}
  \label{figure1}
\end{figure}

We examined numerical performances of our procedure based on 100 Monte Carlo simulations.
All of the tuning parameters $\lambda_{\beta}$, $\lambda_{\gamma}$, $\lambda_{\theta}$ and $\tau_{\theta}$
were selected by the Bayesian information criteria (BIC; \citet{schwarz78}).
For instance, let us consider the selection of $(\lambda_{\beta}, \lambda_{\gamma})$ by BIC.
Let $\hat{\bm{\beta}}(\lambda_{\beta},\lambda_{\gamma})$ and $\hat{\bm{\gamma}}(\lambda_{\beta},\lambda_{\gamma})$
be the outputs of Algorithm \ref{algorithm1} with the tuning parameters $(\lambda_{\beta},\lambda_{\gamma})$.
Then, the optimal tuning parameters $(\hat{\lambda}_{\beta}, \hat{\lambda}_{\gamma})$ are given by
\begin{align*}
(\hat{\lambda}_{\beta}, \hat{\lambda}_{\gamma}) = \argmin_{\lambda_{\beta}>0;\lambda_{\gamma}>0}&
\frac{1}{2n}\|\bm{y} - \bm{X}\hat{\bm{\beta}}(\lambda_{\beta},\lambda_{\gamma}) - \sqrt{n}\hat{\bm{\gamma}}(\lambda_{\beta},\lambda_{\gamma})\|_{2}^{2} \\
&\qquad + \frac{\log n}{n}\big\{|{\rm supp}(\hat{\bm{\beta}}(\lambda_{\beta},\lambda_{\gamma}))| + 
|{\rm supp}(\hat{\bm{\gamma}}(\lambda_{\beta},\lambda_{\gamma}))|\big\}.
\end{align*}
Practically, since it is impossible to search all possible tuning parameters, 
we searched them over some candidate values which were generated by the similar way to
in \citet{friedman10}.

The first two simulations were designed to see the impact of the number of outliers.
The two scenarios $(n,p,s^{*}) = (200,200,10)$ (moderate dimension) and
$(n,p,s^{*}) = (200,400,20)$ (high dimension) with various $g^{*}$ were considered.
The covariates $\bm{x}_{i}$'s were independently drawn from $N_{p}(\bm{0},\bm{\Sigma})$ with
$(\bm{\Sigma})_{ij} = 0.3^{|i-j|}$, the true coefficients were given by $\beta_{j}^{*} = {\rm sgn}(u_{j})$ and
the true outliers were given by $\sqrt{n}\gamma_{i}^{*} = 8$.
The positions of the non-zero coefficients and outliers 
were uniformly drawn from $\{1,\dots, p\}$ and $\{1,\dots, n\}$, respectively.
The $\varepsilon_{i}$'s and $u_{j}$'s were independently drawn from $N(0,1)$.
Moreover, we used $R_{w}=100$ in (\ref{rr}) and Algorithm \ref{algorithm1} stopped when
$\|\bm{\beta}^{k} - \bm{\beta}^{k-1}\|_{1}/\tilde{s} \le 10^{-3}$ was satisfied
at the iteration $k$.

\begin{table}[ht]
\begin{center}
\footnotesize
 \begin{tabular}{cccccccccc}
 \hline
 Pre. & Outlier & Criterion  & Soft & Hard & SCAD &  Garrote & Lasso & Oracle \\ \hline
 (\ref{pre}) & 5\% & $\ell_{2}$-error  & 0.1036 & 0.1008 & 0.1026  & 0.1021 & 1.7400 & 0.3514 \\
  & & FP  & 1.08 & 0.88 & 0.97 & 0.95 & 66.86 & 5.17 \\ 
  & & TP  & 10.00 & 10.00 & 10.00 & 10.00 & 10.00 & 10.00 \\ 
  & 10\% & $\ell_{2}$-error  & 0.1059 & 0.1058 &  0.0992 & 0.1014 & 9.4954 & - \\
  & & FP  & 1.27 & 1.09 & 1.07 & 1.11 & 124.53 & - \\ 
  & & TP  & 10.00 & 10.00 & 10.00 & 10.00 & 9.97 & - \\ 
  & 20\% & $\ell_{2}$-error &0.1828 & 0.1436 & 0.1377 & 0.1485 & 32.0149 &    -  \\ 
  &  & FP & 3.26 & 2.38 & 2.62 & 2.74 & 153.42 &    -  \\
  &  & TP & 10.00 & 10.00 & 10.00 & 10.00 & 9.76 &    -  \\  
  & 30\% & $\ell_{2}$-error  & 0.6788 & 0.4886 & 0.4683 & 0.5107 & 52.2651 & - \\  
  & & FP  & 6.77 & 5.04 & 5.33 & 5.68 & 157.06 & - \\
  & & TP  & 9.87 & 9.89 & 9.89 & 9.89 & 9.53 & - \\  \hline
  (\ref{thpre}) & 5\% & $\ell_{2}$-error & 0.1067 & 0.1060 & 0.1059 & 0.1061 & - & - \\
  & & FP & 1.05 & 1.07 & 1.07 & 1.07 & - & - \\ 
  & & TP & 10.00 & 10.00 & 10.00 & 10.00 & - & - \\ 
  & 10\% & $\ell_{2}$-error & 0.1097 & 0.1080 & 0.1084 &  0.1083  & - & - \\
  & & FP & 1.18 & 1.09 & 1.10 & 1.10  & - & - \\ 
  & & TP & 10.00 & 10.00 & 10.00 & 10.00 & - & - \\ 
  & 20\% & $\ell_{2}$-error & 0.1654 & 0.1501 & 0.1487 & 0.1541 & - &    -  \\ 
  &  & FP & 1.14 & 1.11 & 1.06 & 1.11 & - &    -  \\
  &  & TP & 10.00 & 10.00 & 10.00 & 10.00 & - &    -  \\  
  & 30\% & $\ell_{2}$-error & 0.6583 & 0.5377 & 0.5502 & 0.5788 & - & - \\  
  & & FP & 3.66 & 2.91 & 2.87 & 3.20 & - & - \\
  & & TP & 9.84 & 9.84 & 9.84 & 9.84 & - & - \\   \hline
  \end{tabular}
 \end{center}
 \caption{Numerical performances of the proposed procedure for various outlier percentages 
 when $(n,p,s^{*}) = (200,200,10)$. 
 Each value shows the mean based on 100 Monte Carlo simulations.}
 \label{table1}
 \end{table}

Figure \ref{figure1} shows the support size $\tilde{s} + \tilde{g}$
(left) and the coverage probability $\prob(S^{*}\subset \tilde{S}, G^{*}\subset\tilde{G})$ (right) for the two
preliminary estimators (\ref{pre}) and (\ref{thpre}) when the percentage of outliers increases from 1\% to 35\%.
It can be seen that the two preliminary estimators performed well if the percentage of outliers
was lower than around 20\%. 
The threshold version (\ref{thpre}) had a smaller support size than (\ref{pre}), but
its coverage probability was worse.
We also notice that the coverage probability tended to be low as the percentage of outliers increased.
It would come from the violation of the condition (\ref{screening}).
As the number of outliers increases, the order $a_{n,1}$ increases.

\begin{table}[hb]
\begin{center}
\footnotesize
 \begin{tabular}{cccccccccc}
 \hline
 Pre. & Outlier & Criterion  & Soft & Hard & SCAD &  Garrote & Lasso & Oracle \\ \hline
 (\ref{pre}) & 5\% & $\ell_{2}$-error  & 0.2287 & 0.2139 & 0.2041  & 0.2087 & 3.8340 & 1.0337 \\
  & & FP  & 2.32 & 2.05 & 2.04 & 2.09 & 97.95 & 21.69 \\ 
  & & TP  & 20.00 & 20.00 & 20.00 & 20.00 & 19.98 & 20.00 \\ 
  & 10\% & $\ell_{2}$-error  & 0.3987 & 0.3096 &  0.3039 & 0.3218 & 8.5953 & - \\
  & & FP  & 6.03 & 5.02 & 5.33 & 5.30 & 134.80 & - \\ 
  & & TP  & 20.00 & 20.00 & 20.00 & 20.00 & 19.87 & - \\ 
  & 20\% & $\ell_{2}$-error &3.9041 & 3.5611 & 3.5540 & 3.6990 & 19.5640 &    -  \\ 
  &  & FP & 21.51 & 18.74 & 19.77 & 20.44 & 158.11 &    -  \\
  &  & TP & 19.01 & 18.99 & 19.00 & 19.09 & 18.79 &    -  \\  
  & 30\% & $\ell_{2}$-error &13.7451 & 13.0448 & 13.7788 & 13.9348 & 30.6065 &    -  \\ 
  &  & FP & 43.13 & 40.86 & 42.78 & 43.41 & 169.16 &    -  \\
  &  & TP & 15.85 & 15.81 & 15.90 & 15.82 & 17.52 &    -  \\    \hline
  (\ref{thpre}) & 5\% & $\ell_{2}$-error & 0.2502 & 0.2280 & 0.2232 & 0.2293 & - & - \\
  & & FP & 2.17 & 2.06 & 1.98 & 2.10 & - & - \\ 
  & & TP & 20.00 & 20.00 & 20.00 & 20.00 & - & - \\ 
  & 10\% & $\ell_{2}$-error & 0.4013 & 0.3253 & 0.3137 &  0.3345  & - & - \\
  & & FP & 2.91 & 2.43 & 2.57 & 2.56  & - & - \\ 
  & & TP & 19.99 & 19.99 & 19.99  & 19.99 & - & - \\ 
  & 20\% & $\ell_{2}$-error & 4.1649 & 4.0372 & 3.9793 & 4.0969 & - &    -  \\ 
  &  & FP & 19.78 & 19.04 & 18.81 & 19.04 & - &    -  \\
  &  & TP & 18.99 & 19.01  & 19.00 & 19.00 & - &    -  \\ 
  & 30\% & $\ell_{2}$-error & 13.2793 & 13.9873 & 13.7711 & 14.3158 & - &    -  \\ 
  &  & FP & 41.10 & 41.02 & 41.00 & 41.39 & - &    -  \\
  &  & TP & 16.02 & 16.03  & 15.99 & 16.00 & - &    -  \\  \hline
  \end{tabular}
 \end{center}
 \caption{Numerical performances of the proposed procedure for various outlier percentages 
 when $(n,p,s^{*}) = (200,400,20)$.}
 \label{table2}
 \end{table}

Tables \ref{table1} and \ref{table2} show the squared $\ell_{2}$-error 
$\|\hat{\bm{\beta}} - \bm{\beta}^{*}\|_{2}^{2}$, the number of the false positives
$|\{j|\beta_{j}^{*}=0, \hat{\beta}_{j}\neq 0\}|$ (FP)
and the number of the true positives $|\{j| \beta_{j}^{*}\neq 0, \hat{\beta}_{j}\neq 0\}|$ (TP)
for various penalties for outlier parameters. 
We used the preliminary estimator for $\bm{\beta}^{init}$ in Algorithm \ref{algorithm1}.
We considered the ``Soft", ``Hard", ``SCAD", and ``Garotte" thresholding functions.
Only the Soft does not have the redescending property.
The Garrote has a different behavior from the Hard and the SCAD, in fact,
its $\psi(z;\lambda)$ function never vanishes if $z$ is finite.
For the comparison, we also investigated the performances of the standard ``Lasso" and its ``Oracle" version where
the true outliers are excluded in advance. The symbol ``-" means that
the standard Lasso does not depend on preliminary estimators and its oracle version does not so on
outliers additionally.

As seen in Table \ref{table1}, for the moderate dimension, our procedure provided
quite good estimates and recovered the true support well.
Interestingly, the performance was better than the Oracle.
This would be because the error $\bm{\varepsilon}$ can yield extreme values in the simulation and
the Oracle is not robust against these values.
As seen in Table \ref{table2}, however, 
for the high dimension and large outlier percentages, our procedure did not exhibit good performances.
This result would come from the bad coverage probability of the preliminary estimators.
Compared between the preliminary estimators, (\ref{thpre}) performed better than (\ref{pre}) for the true support recovery,
but the opposite was true for the $\ell_{2}$-error.
It is also noted that the Soft performed worse than the other thresholding functions.
This would be explained by the redescending property.

\begin{figure}[h]
  \centering
  \includegraphics[width=1\textwidth, bb=0 0 799 432]{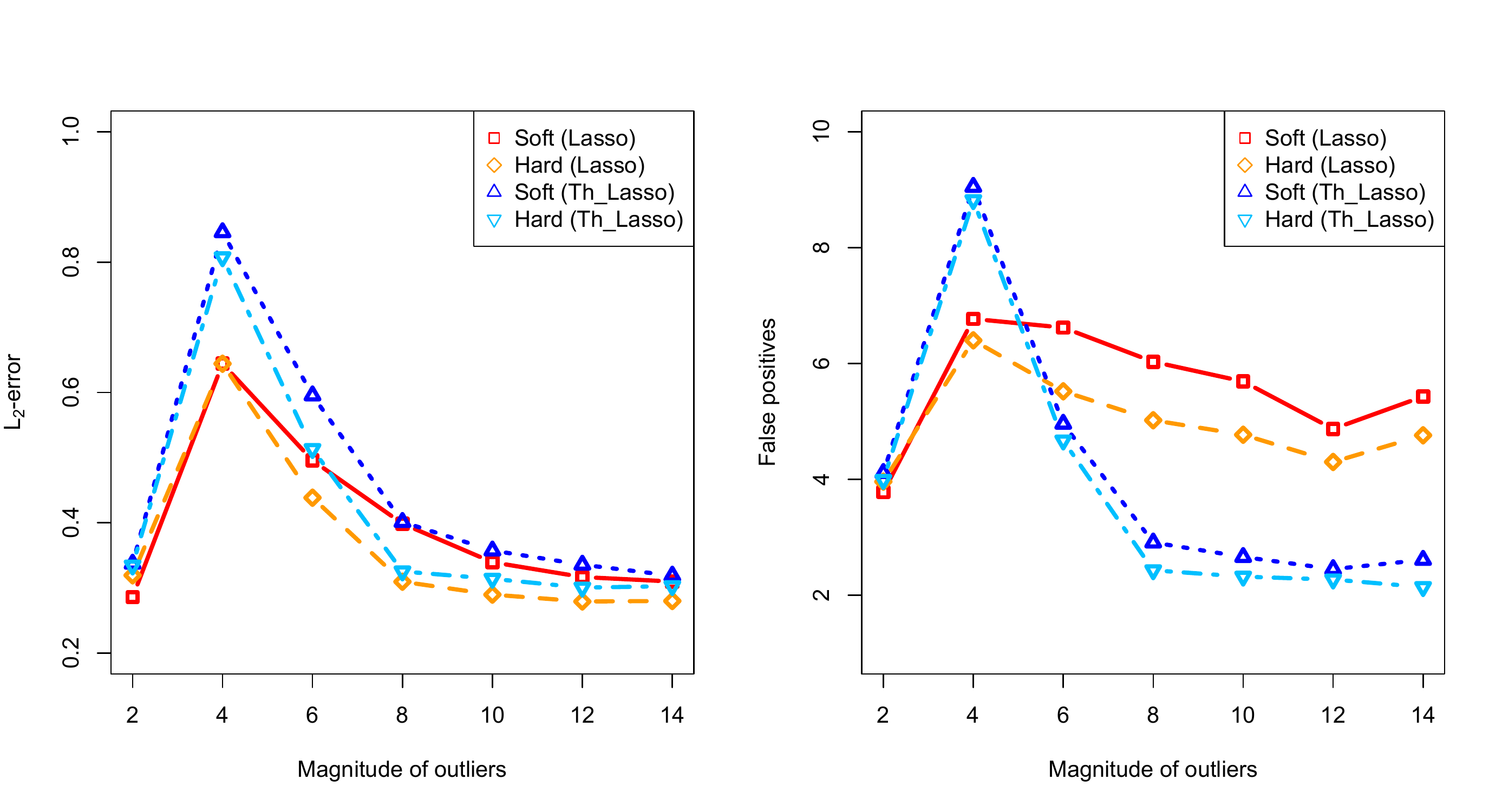}
  \caption{The squared $\ell_{2}$-error (left) and the false positives (right) for various magnitudes of outliers.
  Each point on the curve shows the mean based on 100 Monte Carlo simulations.}
  \label{figure2}
\end{figure}

The final simulations are designed to investigate the impact of the magnitude of outliers.
We used  $(n,p,s^{*}) = (200,400,20)$ with $g^{*}=20$ (10\% outliers) 
and various magnitudes of $\sqrt{n}\bm{\gamma}^{*}$.
We considered the situations $\sqrt{n}\bm{\gamma}^{*}=\gamma^{*}\bm{1}_{n}$ with $\gamma^{*}
\in\{2,4,6,\dots, 14\}$.
The Monte Carlo samples were generated by the same way as above.
In Figure \ref{figure2}, only the performances of the Soft and Hard are shown.
The SCAD and Garotte performed similarly to the Hard.
The true positives are also omitted since they were around 20 for all the cases considered.
Our procedure performed well with low and high magnitudes, but it did not so with a moderate magnitude.
This also would be come from the violation of the condition (\ref{screening}).
For a low magnitude, the outliers would be hidden by the random errors.
When $\bm{\varepsilon}$ is drawn from $N_{n}(\bm{0},\bm{I}_{n})$, the maximum magnitude of $\varepsilon_{i}$'s
is less than $\sqrt{2\log n}$ (it is around 3.3 in this case).
We also note that the performances tended to be stable as the magnitude increased.

\vskip 14pt
\setcounter{chapter}{5}
\setcounter{equation}{0} 
\noindent {\bf 5. Proofs}

\vskip 14pt
\noindent {\bf 5.1. Proof of Proposition \ref{prop2.1}}

Let $(\hat{\bm{\beta}}, \hat{\bm{\gamma}})$ be the output of Algorithm \ref{algorithm1}.
Then, from (\ref{2.3}), we have for $j=1,\dots, p,$
\begin{align*}
\frac{1}{n}\sum_{i=1}^{n}x_{ij}(y_{i}-\bm{x}_{i}^{T}\hat{\bm{\beta}} - \sqrt{n}\hat{\gamma}_{i})
+\lambda_{\beta,j}\partial|\hat{\beta}_{j}| = 0.
\end{align*}
Since $\psi(z;\lambda) = z-\Theta(z;\lambda)$, it follows from (\ref{2.4}) that
\begin{align*}
\frac{1}{n}\sum_{i=1}^{n}x_{ij}(y_{i}-\bm{x}_{i}^{T}\hat{\bm{\beta}} - \sqrt{n}\hat{\gamma}_{i})
&=\frac{1}{n}\sum_{i=1}^{n}x_{ij}(y_{i}-\bm{x}_{i}^{T}\hat{\bm{\beta}} - \Theta(y_{i}-\bm{x}_{i}^{T}\hat{\bm{\beta}};\lambda_{\gamma}w_{\gamma,i})) \\
&=\frac{1}{n}\sum_{i=1}^{n}x_{ij}\psi(y_{i}-\bm{x}_{i}^{T}\hat{\bm{\beta}};\lambda_{\gamma}w_{\gamma,i}),
\end{align*}
which ends the proof.

\vskip 14pt
\noindent {\bf 5.2. Proof of Theorem \ref{th3.1}}

To prove the theorem we prepare the next lemma.
\begin{lemma}\label{le6.1}
Let $Z_{1},\dots, Z_{n}$ be independently identically distributed as the zero mean sub-Gaussian distribution
with a parameter $\sigma > 0$. Then, for any vector $\bm{a}\in\mathbb{R}^{n}$ and any $t \ge 0$,
\begin{align*}
\prob\bigg(\bigg|\sum_{i=1}^{n}a_{i}Z_{i}\bigg| > t\bigg) \le 2\exp\bigg(-\frac{t^{2}}{2\sigma^{2}\|\bm{a}\|_{2}^{2}}\bigg).
\end{align*}
\end{lemma}

\noindent
{\bf Proof.}
Let $Z=\sum_{i=1}^{n}a_{i}Z_{i}$.
From the Markov's inequality, we have for $u >0$
\begin{align*}
\prob\bigg(Z > t\bigg) 
= \prob(e^{uZ} > e^{ut})
\le e^{-ut}\prod_{i=1}^{n}\E(e^{ua_{i}Z_{i}}) 
\le \exp\bigg(-ut + \frac{u^{2}\sigma^{2}\|\bm{a}\|_{2}^{2}}{2}\bigg).
\end{align*}
Note that $u=t/(\sigma^{2}\|\bm{a}\|_{2}^{2})$ minimizes 
$-ut + u^{2}\sigma^{2}\|\bm{a}\|_{2}^{2}/2$ over $u > 0$. 
Thus, $\prob(Z > t)\le \exp\{-t^{2}/(2\sigma^{2}\|\bm{a}\|_{2}^{2})\}$.
Similarly, $\prob(Z <  - t)\le \exp\{-t^{2}/(2\sigma^{2}\|\bm{a}\|_{2}^{2})\}$
and then from $P(|Z| > t)\le P(Z > t) + P(Z < -t)$ we obtain the lemma. 
\qed

\mbox{}

Since we consider the adaptive Lasso type estimator
and the screening property that $S^{*}\subset \tilde{S}$ and $G^{*}\subset \tilde{G}$ is satisfied under
Condition \ref{condition3} and (\ref{screening}),
it suffices to focus only on the covariates selected by the preliminary estimator $\tilde{\bm{\beta}}$, that is,
on the sub-matrix $\bm{X}_{\tilde{S}}=\{x_{ij}|\;1\le i\le n; j \in \tilde{S}\}$.
Correspondingly, we only focus on the coefficients in the set $\tilde{S}$.
In this section we omit the subscript $\tilde{S}$ for simplicity, therefore keep in mind that 
all of the following $\bm{X}$, $\bm{\beta}^{k}$ and $\bm{\beta}^{*}$ have the $\tilde{s}$ dimension.

We shall show the bound (\ref{3.3})
on the event that Condition \ref{condition3} is satisfied and on
\begin{align}\label{6.1}
E=\bigg\{\bigg\|\frac{1}{n}\bm{X}^{T}_{(\tilde{G}^{c})}\bm{\varepsilon}_{\tilde{G}^{c}}\bigg\|_{\infty}\le C\sqrt{\frac{\sigma^{2}\log p}{n}}\bigg\}
\end{align}
both of which have probabilities going to one, 
where $\bm{X}_{(\tilde{G}^{c})} = \{x_{ij}|\; i\in \tilde{G}^{c}, j\in \tilde{S}\}$,
$\e_{\tilde{G}^{c}} = \{\varepsilon_{i}|\;i\in\tilde{G}^{c}\}$ and $C>\sqrt{2}$.
In fact, from Lemma \ref{le6.1} and Condition \ref{condition1}, we have
for given the preliminary estimators $\tilde{\bm{\beta}}$ and $\tilde{\bm{\gamma}}$,
\begin{align*}
\prob(E) &= 1- \prob(E^{c})=1-\prob\bigg\{\bigg\|\frac{1}{n}\bm{X}^{T}_{(\tilde{G}^{c})}\bm{\varepsilon}_{\tilde{G}^{c}}\bigg\|_{\infty}> C\sqrt{\frac{\sigma^{2}\log p}{n}}\bigg\} \\
&\ge 1- \sum_{j \in \tilde{S}}\prob\bigg(\bigg|\frac{1}{n}\sum_{i\in\tilde{G}^{c}}x_{ij}\varepsilon_{i}\bigg| > C\sqrt{\frac{\sigma^{2}\log p}{n}}\bigg) \\
&\ge 1- 2\sum_{j\in\tilde{S}}\exp(-2^{-1}C^{2}\log p) \\
&\ge 1- 2\exp(\log p - 2^{-1}C^{2}\log p) = 1-o(1).
\end{align*}
Note that the lower bound does not depend on the preliminary estimators, and hence 
the probability of $E$ goes to one.

Since $\bm{\beta}^{k} = {\rm argmin}_{\bm{\beta}}L(\bm{\beta},\bm{\gamma}^{k-1})$, 
it follows that $L(\bm{\beta}^{k},\bm{\gamma}^{k-1}) \le L(\bm{\beta}^{*}, \bm{\gamma}^{k-1})$.
Hence, we have 
\begin{align*}
\frac{1}{2n}\|\bm{X}\bm{\Delta}^{k}\|_{2}^{2}
&\le \frac{1}{n}\langle\bm{X}^{T}(\bm{y}-\bm{X}\bm{\beta}^{*} - \sqrt{n}\bm{\gamma}^{k-1}), \bm{\Delta}^{k}\rangle
+\lambda_{\beta}\sum_{j\in\tilde{S}}w_{\beta,j}(|\beta_{j}^{*}| - |\beta_{j}^{k}|) \\
&=A_{1} + A_{2},\quad {\rm say}, 
\end{align*}
where $\bm{\Delta}^{k} = \bm{\beta}^{k} - \bm{\beta}^{*}$.
First, we evaluate the term $A_{1}$.
Since $\sqrt{n}\gamma_{i}^{k-1} = \Theta(y_{i} - \bm{x}_{i}^{T}\bm{\beta}^{k-1};\lambda_{\gamma}w_{\gamma,i})$
and $\sqrt{n}\gamma_{i}^{k-1} = 0$ when $i \in \tilde{G}^{c}$, 
the $j$-th element of $\bm{X}^{T}(\bm{y} - \bm{X}\bm{\beta}^{*}-
\sqrt{n}\bm{\gamma}^{k-1})$ is given by
\begin{align*}
&\sum_{i\in\tilde{G}}x_{ij}\{y_{i} - \bm{x}_{i}^{T}\bm{\beta}^{*}-\Theta(y_{i}-\bm{x}_{i}^{T}\bm{\beta}^{k-1};\lambda_{\gamma}w_{\gamma,i})\} + \sum_{i\in\tilde{G}^{c}}x_{ij}(y_{i} - \bm{x}_{i}^{T}\bm{\beta}^{*}) \\
&\quad =
\sum_{i\in\tilde{G}}x_{ij}\big\{\bm{x}_{i}^{T}\bm{\Delta}^{k-1} + (y_{i}-\bm{x}_{i}^{T}\bm{\beta}^{k-1})
-\Theta(y_{i}-\bm{x}_{i}^{T}\bm{\beta}^{k-1};\lambda_{\gamma}w_{\gamma,i})\big\}  \\
&\qquad \qquad + \sum_{i\in\tilde{G}^{c}}x_{ij}(y_{i} - \bm{x}_{i}^{T}\bm{\beta}^{*}).
\end{align*}
Thus, we obtain 
\begin{align*}
A_{1}&=\frac{1}{n}(\bm{X}_{(\tilde{G})}\bm{\Delta}^{k-1})^{T}(\bm{X}_{(\tilde{G})}\bm{\Delta}^{k}) \\
&\quad
+\frac{1}{n}\sum_{j\in\tilde{S}}\Delta_{j}^{k}\sum_{i\in\tilde{G}}x_{ij}\{(y_{i}-\bm{x}_{i}^{T}\bm{\beta}^{k-1})
-\Theta(y_{i}-\bm{x}_{i}^{T}\bm{\beta}^{k-1};\lambda_{\gamma}w_{\gamma,i})\} \\
&\quad
+\frac{1}{n}\sum_{j\in\tilde{S}}\Delta_{j}^{k}\sum_{i\in\tilde{G}^{c}}x_{ij}(y_{i}-\bm{x}_{i}^{T}\bm{\beta}^{*})
=A_{11}+A_{12}+A_{13},\quad{\rm say}.
\end{align*}
Since $\|\bm{\Delta}^{k}\|_{0}\le\tilde{s}$ for any $k\ge 1$, the definition of
the (doubly) restricted largest eigenvalue in (\ref{3.2}) implies that
\begin{align}\label{6.2}
|A_{11}| \le \frac{1}{n}\|\bm{X}_{(\tilde{G})}\bm{\Delta}^{k-1}\|_{2}\|\bm{X}_{(\tilde{G})}\bm{\Delta}^{k}\|_{2}
\le \delta_{max}(\tilde{s},\tilde{g})\|\bm{\Delta}^{k-1}\|_{2}\|\bm{\Delta}^{k}\|_{2}.
\end{align}
Note that $|\Theta(x;\lambda) - x| \le \lambda$ under Condition 
\ref{condition2}. Then,
\begin{align}\label{6.3}
|A_{12}| \le \frac{R_{w}\lambda_{\gamma}}{n}\|\bm{X}_{(\tilde{G})}\|_{\ell_{1}}\|\bm{\Delta}^{k}\|_{1}
\le C\sqrt{\frac{\sigma^{2}\log p}{n}}\|\bm{\Delta}^{k}\|_{1},
\end{align}
where $R_{w} > 0$ is defined in (\ref{rr}).
Since $G^{*} \subset \tilde{G}$, we have $\gamma_{i}^{*} = 0$ for $i\in\tilde{G}^{c}$. Then, 
from (\ref{6.1}), we have
\begin{align}\label{6.4}
|A_{13}|\le \frac{1}{n}\|\bm{X}_{(\tilde{G}^{c})}\bm{\varepsilon}_{\tilde{G}^{c}}\|_{\infty}\|\bm{\Delta}^{k}\|_{1}
\le C\sqrt{\frac{\sigma^{2}\log p}{n}}\|\bm{\Delta}^{k}\|_{1}.
\end{align}

For the term $A_{2}$, since $S^{*}\subset \tilde{S}$, 
\begin{align*}
A_{2}
&\le\lambda_{\beta}\sum_{j\in S^{*}}w_{\beta,j}|\beta_{j}^{*}|
-\lambda_{\beta}\sum_{j\in S^{*}}w_{\beta,j}|\beta_{j}^{*} + \Delta_{j}^{k}|
-\lambda_{\beta}\sum_{j\in S^{*^{c}}}w_{\beta,j}|\beta_{j}^{*} + \Delta_{j}^{k}| \\
&\le
 \lambda_{\beta}\sum_{j\in S^{*}}w_{\beta,j}|\Delta_{j}^{k}|
 -\lambda_{\beta}\sum_{j\in S^{*^{c}}}w_{\beta,j}|\Delta_{j}^{k}| \\
 &\le \lambda_{\beta}\max_{j\in S^{*}}w_{\beta,j}\|\bm{\Delta}_{S^{*}}^{k}\|_{1}
 -\frac{\lambda_{\beta}}{R_{w}}\|\bm{\Delta}_{S^{*^{c}}}^{k}\|_{1}
\end{align*}
Combined with (\ref{6.2})--(\ref{6.4}), it follows from 
$\|\bm{\Delta}^{k}\|_{1}=\|\bm{\Delta}_{S^{*}}^{k}\|_{1} + \|\bm{\Delta}_{S^{*c}}^{k}\|_{1}$ that
\begin{align*}
\frac{1}{2n}\|\bm{X}\bm{\Delta}^{k}\|_{2}^{2}
&\le \delta_{max}(\tilde{s},\tilde{g})\|\bm{\Delta}^{k-1}\|_{2}\|\bm{\Delta}^{k}\|_{2}
+2C\sqrt{\frac{\sigma^{2}\log p}{n}}\|\bm{\Delta}^{k}\|_{1} \\
&\qquad + \lambda_{\beta}\max_{j\in S^{*}}w_{\beta,j}\|\bm{\Delta}_{S^{*}}^{k}\|_{1}
 -\frac{\lambda_{\beta}}{R_{w}}\|\bm{\Delta}_{S^{*c}}^{k}\|_{1} \\
 &=\delta_{max}(\tilde{s},\tilde{g})\|\bm{\Delta}^{k-1}\|_{2}\|\bm{\Delta}^{k}\|_{2}
 +\left(\lambda_{\beta}\max_{j\in S^{*}}w_{\beta,j} + 2C\sqrt{\frac{\sigma^{2}\log p}{n}}\right)
 \|\bm{\Delta}^{k}_{S^{*}}\|_{1} \\
 &\qquad + \left(2C\sqrt{\frac{\sigma^{2}\log p}{n}} - \frac{\lambda_{\beta}}{R_{w}}\right)\|\bm{\Delta}_{S^{*c}}^{k}\|_{1} \\
 &\qquad \le \delta_{max}(\tilde{s},\tilde{g})\|\bm{\Delta}^{k-1}\|_{2}\|\bm{\Delta}^{k}\|_{2} +
 \lambda_{\beta}\big(R_{w}^{-1}+\max_{j\in S^{*}}w_{\beta,j}\big)\sqrt{s^{*}}\|\bm{\Delta}^{k}\|_{2}.
\end{align*}
From Condition \ref{condition3} and (\ref{3.1}), we have 
$\frac{1}{n}\|\bm{X}\bm{\Delta}^{k}\|_{2}^{2} \ge \kappa^{}\|\bm{\Delta}^{k}\|_{2}^{2}$. Thus, for $k\ge 1$,
\begin{align*}
\|\bm{\Delta}^{k}\|_{2}^{} \le \rho\|\bm{\Delta}^{k-1}\|_{2} + 2\kappa^{-1}\sqrt{s^{*}}\lambda_{\beta}\big(R_{w}^{-1}+\max_{j\in S^{*}}w_{\beta,j}\big).
\end{align*}
Let $\bm{\Delta}^{0} = \bm{\beta}^{init} - \bm{\beta}^{*}$. Then, the bound (\ref{3.3}) is derived by
solving the above recurrence relation for $k=1,2,\dots$, which ends the proof.

\vskip 14pt
\noindent {\bf 5.3. Proof of Theorem \ref{th3.2}}

Throughout this section we denote positive constants by $C_{i}\;(i\ge 1)$ which may be different from each other.
Suppose that for some $\ell \in \tilde{S}\cap S^{*^{c}}$, $\beta_{\ell}^{k}\neq 0$. 
Without loss of generality, we can assume $\beta_{\ell}^{k} > 0$.
Then, from the first order condition for $\bm{\beta}^{k}$, the value
\begin{align}\label{6.5}
\frac{1}{n}\sum_{j\in\tilde{S}}\sum_{i=1}^{n}x_{ij}x_{i\ell}\beta_{j}^{k}
-\frac{1}{n}\sum_{i=1}^{n}x_{i\ell}\bigg\{y_{i}-\Theta\big(y_{i} - \sum_{j\in\tilde{S}}x_{ij}\beta_{j}^{k-1};\lambda_{\gamma}w_{\gamma,i}\big)\bigg\}+ \lambda_{\beta} w_{\beta,\ell}
\end{align}
should be zero.
But, if we can show the first two terms are dominated by the third term $\lambda_{\beta}w_{\beta,\ell}$
for any $\ell \in  \tilde{S}\cap S^{*^{c}}$, 
then the above value can never be zero, which leads to the contradiction.

First, we evaluate the middle term of (\ref{6.5}). From the definition, the inside of $\Theta$ is represented as
\begin{align*}
y_{i} - \sum_{j\in\tilde{S}}x_{ij}\beta_{j}^{k-1}=
\begin{cases}
\sqrt{n}\gamma_{i}^{*} + \varepsilon_{i} -\sum_{j\in\tilde{S}}x_{ij}(\beta_{j}^{k-1} - \beta_{j}^{*}), & i\in G^{*} \cap \tilde{G} = G^{*} \\
 \varepsilon_{i} -\sum_{j\in\tilde{S}}x_{ij}(\beta_{j}^{k-1} - \beta_{j}^{*}), & i\in G^{*^{c}} \cap \tilde{G}.
\end{cases}
\end{align*}
By Lemma \ref{le6.1}, we can show that $\max_{1\le i\le n}|\varepsilon_{i}| \le C_{1}\sqrt{\log n}$
with probability going to one.
Let $\lambda_{\beta}=C_{2}\{(\log p)/n\}^{1/2}$,
then it follows from Corollary \ref{cor3.1} and $a_{n,2}s^{*} = o(n)$ that
\begin{align*}
\max_{1\le i\le n}\bigg|\varepsilon_{i} -\sum_{j\in\tilde{S}}x_{ij}(\beta_{j}^{k-1} - \beta_{j}^{*})\bigg|
&\le \max_{1\le i\le n}|\varepsilon_{i}| + C_{3}\sqrt{\tilde{s}}\|\bm{\beta}^{k-1} - \bm{\beta}^{*}\|_{2} \\
&\le C_{4}\sqrt{\log n}(1 + o(1)),
\end{align*}
which implies that if 
$\lambda_{\gamma}\min_{i\in G^{*^{c}} \cap \tilde{G}}w_{\gamma, i} \ge C_{4}\sqrt{\log n}$, then
$y_{i}-\sum_{j\in\tilde{S}}x_{ij}\beta_{j}^{k-1}=0$ for $i \in G^{*^c}\cap \tilde{G}$. 
Note that $\min_{i\in G^{*^{c}} \cap \tilde{G}}w_{\gamma, i} = R_{w}$ for sufficiently large $n$ since 
$\min_{i\in G^{*^c}\cap \tilde{G}}(1/|\tilde{\gamma}_{i}|) \ge 1/(Ca_{n,1}) \to \infty$.
Hence it suffices to put $\lambda_{\gamma} = C_{5}\sqrt{\log n}$ for a sufficiently large $C_{5} > 0$.
Such a $\lambda_{\gamma}$ satisfying the condition of Corollary \ref{cor3.1} can be selected since $a_{n,2}^{2}=o(n)$.
Meanwhile, since $\sqrt{\log n}=o(\sqrt{n}\min_{i\in G^{*}}|\gamma_{i}^{*}|)$, we have
$y_{i}-\sum_{j\in\tilde{S}}x_{ij}\beta_{j}^{k-1}=\sqrt{n}\gamma_{i}^{*}(1+o(1))$ for $i\in G^{*}$. 
Thus, for sufficiently large $n$, it holds that 
$\min_{i\in G^{*}}|y_{i} - \sum_{j\in\tilde{S}}x_{ij}\beta_{j}^{k-1}| \ge \lambda_{\gamma} R_{w} \ge \lambda_{\gamma}\max_{i\in G^{*}}w_{\gamma,i}$.
Therefore, under Condition \ref{condition2},
\begin{align*}
\Theta\big(y_{i} - \sum_{j\in\tilde{S}}x_{ij}\beta_{j}^{k-1};\lambda_{\gamma}w_{\gamma,i}\big)
=\begin{cases}
y_{i} - \sum_{j\in\tilde{S}}x_{ij}\beta_{j}^{k-1} + O(\sqrt{\log n}), & i\in G^{*} \\
0, & i\in G^{*^c}\cap \tilde{G}.
\end{cases}
\end{align*}
Then, it follows that
\begin{align*}
&\sum_{i=1}^{n}x_{i\ell}\bigg\{y_{i} - \Theta\big(y_{i}-\sum_{j\in\tilde{S}}x_{ij}\beta_{j}^{k-1};\lambda_{\gamma}w_{\gamma,i}\big)\bigg\} \\
&\qquad =
\sum_{i\in \tilde{G}}x_{i\ell}\bigg\{y_{i} - \Theta\big(y_{i}-\sum_{j\in\tilde{S}}x_{ij}\beta_{j}^{k-1};\lambda_{\gamma}w_{\gamma,i}\big)\bigg\} + \sum_{i\in\tilde{G}^{c}}x_{i\ell}y_{i} \\
&\qquad =
\sum_{i\in G^{*}}x_{i\ell}\sum_{j\in\tilde{S}}x_{ij}\beta_{j}^{k-1} + O(\sqrt{\log n})\sum_{i\in G^{*}}x_{i\ell}
+ \sum_{i\in G^{*^c}\cap \tilde{G}}x_{i\ell}y_{i}+ \sum_{i\in\tilde{G}^{c}}x_{i\ell}y_{i} \\
&\qquad =
\sum_{i\in G^{*}}x_{i\ell}\sum_{j\in\tilde{S}}x_{ij}\big\{\beta_{j}^{*} + (\beta_{j}^{k-1} - \beta_{j}^{*})\big\} + \sum_{i\in G^{*^c}}x_{i\ell}y_{i} + O(g^{*}\sqrt{\log n}) \\
&\qquad =
\sum_{i=1}^{n}\sum_{j\in\tilde{S}}x_{i\ell}x_{ij}\beta_{j}^{*} + \sum_{i\in G^{*}}\sum_{j\in\tilde{S}}x_{i\ell}x_{ij}(\beta_{j}^{k-1} - \beta_{j}^{*}) + \sum_{i\in G^{*^c}}x_{i\ell}\varepsilon_{i} + O(g^{*}\sqrt{\log n}).
\end{align*}
Thus, (\ref{6.5}) can be written as 
\begin{align*}
&\frac{1}{n}\sum_{i=1}^{n}\sum_{j\in\tilde{S}}x_{ij}x_{i\ell}(\beta_{j}^{k} - \beta_{j}^{*}) -\frac{1}{n}\sum_{i\in G^{*}}\sum_{j\in\tilde{S}}x_{ij}x_{i\ell}(\beta_{j}^{k-1} - \beta_{j}^{*})  \\
&\qquad\qquad - \frac{1}{n} \sum_{i\in G^{*^c}}x_{i\ell}\varepsilon_{i} +O\bigg(\frac{g^{*}}{n}\sqrt{\log p}\bigg)+ \lambda_{\beta}w_{\beta,\ell}.
\end{align*}
Clearly the fist and second terms have the order $\sqrt{a_{n,2}s^{*}}\lambda_{\beta}$.
From the proof of Theorem \ref{th3.1}, the third term is of order $\{(\log p)/n\}^{1/2}$.
By  $a_{n,1}\max\big(\sqrt{a_{n,2}s^{*}},g^{*}/\sqrt{n}\big) = o(1)$, the first four terms of (\ref{6.5}) are dominated by 
$\lambda_{\beta}\min_{j\in \tilde{S}\cap S^{*^c}}w_{\beta,j}$
since $\min_{j\in \tilde{S}\cap S^{*^c}}w_{\beta,j} \ge 1/(Ca_{n,1})$.



\vskip 14pt
\noindent {\large\bf Acknowledgements}

This work was supported by the System Genetics Project of the
Research Organization of Information and Systems.
\par

\markboth{\hfill{\footnotesize\rm SHOTA KATAYAMA AND HIRONORI FUJISAWA} \hfill}
{\hfill {\footnotesize\rm SPARSE AND ROBUST LINEAR REGRESSION} \hfill}

\bibhang=1.7pc
\bibsep=2pt
\fontsize{9}{14pt plus.8pt minus .6pt}\selectfont
\renewcommand\bibname
\vskip 1.5cm
\noindent
Tokyo Institute of Technology, Tokyo, Japan.
\vskip 2pt
\noindent
E-mail: katayama.s.ad@m.titech.ac.jp
\vskip 2pt

\noindent
The Institute of Statistical Mathematics, Tokyo, Japan.
\vskip 2pt
\noindent
E-mail: fujisawa@ism.ac.jp
\end{document}